\documentclass[12pt]{article}
\usepackage{amsmath}
\usepackage{amsfonts}
\usepackage{amssymb}
\usepackage[mathscr]{eucal}

\newif{\ifcomentarios}
\comentariosfalse

\newtheorem{theorem}{Theorem}

\newtheorem{corollary}[theorem]{Corollary}

\newtheorem{definition}[theorem]{Definition}

\newtheorem{lemma}[theorem]{Lemma}

\newtheorem{remark}[theorem]{Remark}

\newtheorem{proposition}[theorem]{Proposition}

\newcommand{\be}{\begin{eqnarray}}
\newcommand{\en}{\end{eqnarray}}
\newcommand{\bee}{\begin{eqnarray*}}
\newcommand{\ene}{\end{eqnarray*}}

\begin{document}

\title{Sparse Block--Jacobi Matrices with Exact Hausdorff Dimension}
\author{S. L. Carvalho\thanks{
Supported by FAPESP under grant \#06/60711-4. Email: \texttt{
silas@fig.if.usp.br}}, \ D. H. U. Marchetti\thanks{
No thanks. Email: \texttt{marchett@if.usp.br}} \ \& W. F. Wreszinski\thanks{
Partially supported by CNPq. Email:\texttt{\ wreszins@fma.if.usp.br}} \\
Instituto de F\'{\i}sica\\
Universidade de S\~{a}o Paulo \\
Caixa Postal 66318\\
05314-970 S\~{a}o Paulo, SP, Brasil}

\date{}
\maketitle

\begin{abstract}
We show that the Hausdorff dimension of the spectral measure of a class of
deterministic, i. e. nonrandom, block--Jacobi matrices may be determined
exactly, improving a result of Zlato\v{s} (J. Funct. Anal. \textbf{207},
216-252 (2004)).
\end{abstract}


\section{Introduction\label{2}}

\setcounter{equation}{0} \setcounter{theorem}{0}

In \cite{MarWre} two of the authors addressed the possibility that a
spectral transition takes place in a deterministic model. The model is
represented by a class of Jacobi matrices with a sparse potential in the
sense that the perturbation of the free Jacobi matrix (the $0$--Dirichlet
Laplacean on $l^{2}(\mathbb{Z}_{+})$) is a (direct) sum of a fixed $2\times
2 $ off--diagonal matrix placed at sites whose distances from one another
grow exponentially. In the present work we improve and complement results of 
\cite{MarWre} in two directions. The model is now represented by
block--Jacobi matrices and we are able to compute, for sparse perturbations
satisfying transversal homogeneity, the exact Hausdorff dimension of their
spectral measures.

Denoting the set of non--negative integers by $\mathbb{Z}_{+}$, let $\Lambda
=\mathbb{Z}_{+}\times \{0,1,\ldots ,L-1\}$ be a strip of width $L\geq 1$ on
the $\mathbb{Z}_{+}^{2}$ plane and define, on the separable Hilbert space $%
l^{2}(\Lambda ,\mathbb{C})$, an operator $\Delta _{P,\phi }$ for each
sequence $P=(p_{n})_{n\geq -1}$ of numbers $p_{n}\in (0,1]$ and angle $\phi
\in \lbrack 0,\pi )$: 
\begin{equation}
\!\!\!(\Delta _{P,\phi
}u)(k,m):=p_{k}u(k+1,m)+p_{k-1}u(k-1,m)+u(k,m+1)+u(k,m-1)\;,  \label{Pip}
\end{equation}%
\noindent for all $(k,m)\in \Lambda $ with phase boundary conditions at $%
k=-1 $: 
\begin{equation}
u(-1,m)\cos \phi -u(0,m)\sin \phi =0  \label{boucon}
\end{equation}%
\noindent for each $m\in \{0,\ldots ,L-1\}$, and periodic boundary
conditions on the vertical direction: $u(k,L)=u(k,0)$ for each $k\in \mathbb{%
Z}_{+}$.

The operator $\Delta _{P,0}$ with phase boundary $0$ is, in particular,
defined on a cylinder with the $0$--Dirichlet boundary condition on $k=-1$: $%
u(-1,m)=0$ for every $m\in \{0,\ldots ,L-1\}$ and the operator $\Delta
_{0,\phi }$, defined by setting $p_{n}=1$, $\forall n$ , reduces to the
usual discrete Laplacean in $\Lambda $ with phase boundary $\phi $. We note
that $p_{k}$ lives on the horizontal edges \thinspace $\left\langle
(k,m),(k+1,m)\right\rangle $, $m\in \{0,\ldots ,L-1\}$, and (\ref{Pip}) is
defined with the same value $p_{k}$ to each $m$. Such property is referred
to in the present text as \textit{transversal homogeneity}.

The sparse perturbation considered here is a natural extension of the
perturbation employed on the one-dimensional problem developed on \cite%
{MarWre}. By sparse perturbation we mean a perturbation about the Laplacean: 
$\Delta _{P,\phi }=\Delta _{0,\phi }+V_{P}$ where the potential $V_{P}$ is
composed of infinitely many vertical `barriers' whose distances from one
another grow exponentially. The sequence $P=(p_{n})_{n\geq -1}$ of
`barriers' is of the form 
\begin{equation}
p_{n}=\left\{ 
\begin{array}{lll}
1-\delta  & \mathrm{if} & n=a_{j}\in \mathcal{A}\,, \\ 
1 & \mathrm{if} & n\not\in \mathcal{A}\,,%
\end{array}%
\right.   \label{Pe}
\end{equation}%
\noindent for $\delta \in (0,1)$ and a set of positive integers $\mathcal{A}%
=\{a_{j}\}_{j\geq 1}$ such that 
\begin{equation}
a_{j}-a_{j-1}\geq 2,\qquad \qquad j=2,3,\ldots   \label{isolate}
\end{equation}%
\noindent and 
\begin{equation*}
\lim_{j\rightarrow \infty }\frac{a_{j+1}}{a_{j}}=\beta >1\;.
\end{equation*}%
\noindent Condition (\ref{isolate}) makes each `barrier' to be located in an
isolated single column of horizontal edges and $\beta $ is the so called
\textquotedblleft sparseness parameter\textquotedblright . As in \cite%
{MarWre}, the separations between the barriers are fixed as 
\begin{equation}
a_{j}-a_{j-1}=\beta ^{j},\qquad \qquad j=2,3,\ldots   \label{cesp}
\end{equation}%
with $a_{1}=\beta \geq 2$ an integer, in order to simplify our analysis.
From now on, (\ref{Pe}) with $\mathcal{A}$ given by (\ref{cesp}) will be the
only sequence considered and we shall denote by $\Delta _{\delta ,\phi }$
the corresponding operator with $P=(p_{n})_{n\geq -1}$ of this form.

The operator $\Delta _{P,\phi }$ with $\phi =0$ may be written in the
block--Jacobi matrix form 
\begin{equation}
\mathcal{J}_{P}=J_{P}\otimes I_{L}+I\otimes A_{L}\;,  \label{kronecker}
\end{equation}%
\noindent with $I$ the identity operator on $l^{2}(\mathbb{Z}_{+})$, $J_{P}$
is defined by 
\begin{equation}
J_{P}=\left( 
\begin{array}{ccccc}
0 & p_{0} & 0 & 0 & \cdots \\ 
p_{0} & 0 & p_{1} & 0 & \cdots \\ 
0 & p_{1} & 0 & p_{2} & \cdots \\ 
0 & 0 & p_{2} & 0 & \cdots \\ 
\vdots & \vdots & \vdots & \vdots & \ddots%
\end{array}%
\right) \;,  \label{MAJA}
\end{equation}%
\noindent the $(p_{n})_{n\geq 0}$ as in (\ref{Pe}), $A_{L}$ and $I_{L}$
denoting, respectively, the $L\times L$ matrix 
\begin{equation*}
A_{L}=\left( 
\begin{array}{cccccc}
0 & 1 & 0 & \cdots & 0 & 1 \\ 
1 & 0 & 1 & \cdots & 0 & 0 \\ 
\vdots & \vdots & \vdots & \ddots & \vdots & \vdots \\ 
0 & 0 & 0 & \cdots & 0 & 1 \\ 
1 & 0 & 0 & \cdots & 1 & 0%
\end{array}%
\right) \;
\end{equation*}%
and the $L\times L$ identity. The matrix class above reduces to the one
studied on \cite{MarWre} by setting the strip width $L=1$ ($I_{L}=1$ and $%
A_{L}=0$ in this case).

It is interesting to note that as $\delta $ varies in the interval $(0,1)$, $%
J_{\delta }$ (given by $J_{P}$ with $P$ satisfying (\ref{Pe})) interpolates
continuously two distinct situations: a dense pure point spectrum at $\delta
=1$ and an absolutely continuous spectrum at $\delta =0$. For a more
detailed discussion which includes a wider class of perturbations see
Section 1 of \cite{MarWre}. 

\begin{remark}
The results presented in this paper are not restricted to the operator $%
\Delta _{\delta ,\phi }$. We could extend our methods to any sparse
perturbation that is block-diagonalizable, i.e., that can be decomposed into
its one-dimensional constituents by a matrix conjugation. The transversal
homogeneity condition allows us to use the discrete Fourier transform to
reduce to this form. The results also hold if $\phi $ in (\ref{boucon}) is
different for each $m$.
\end{remark}

\begin{remark}
The operator $\Delta _{\delta ,\phi }$ with $\phi $--phase boundary
condition at $k=-1$ (\ref{boucon}) may also be written in the block--Jacobi
matrix form (\ref{kronecker}). If $\mathcal{J}_{\delta ,\phi }$ denotes the
corresponding matrix, we have 
\begin{equation}
\mathcal{J}_{\delta ,\phi }=\mathcal{J}_{\delta }+E_{0}\otimes \tan \phi
I_{L}\;,  \label{JJ+}
\end{equation}%
\noindent where $E_{0}$ is an operator on $l^{2}(\mathbb{Z}_{+})$ with all
elements zero except $\left( E_{0}\right) _{00}=1$. If the $\phi $--phase
condition varies for each $m$, $\tan \phi I_{L}$ in (\ref{JJ+}) is replaced
by $\text{diag}\left\{ \tan \phi _{m}\right\} _{m=0}^{L-1}$.
\end{remark}

The very basic method employed to study the spectrum of sparse Schr\"{o}%
dinger operators is given by Pearson \cite{Pearson}. Let $\Delta _{\delta
,0}^{k}$ be the sparse operator $\Delta _{\delta ,0}$ ($\Delta _{\delta
,0}=J_{\delta }$ for $L=1$) with $(p_{n})_{n}$ given by (\ref{Pe}) if $%
n<a_{k}$ and $p_{n}=1$ for all $n\geq a_{k}$, and let $\rho _{k}(\varphi )$
denote the corresponding spectral measure. Note that $d\rho _{k}/d\varphi $
exists for almost every $\varphi \in \left[ 0,\pi \right] $ and $\rho _{k}$
is absolutely continuous with respect to the Lebesgue measure. The spectral
measure $\rho $ of $\Delta _{\delta ,0}$, which may be derived from the
limit as $k\rightarrow \infty $ of $\rho _{k}$, is determined by the
asymptotic behavior, as $n\rightarrow \infty $, of the solution $\psi
_{n}=\psi _{n}(\varphi )$ of the equation%
\begin{equation*}
\left( \Delta _{\delta ,0}^{k}\psi \right) _{n}=\lambda \psi _{n}~,\qquad
\lambda =2\cos \varphi 
\end{equation*}%
in the following sense. If $R_{k}(\varphi )$ and $\theta _{k}(\varphi )$ are
the radius and angle of Pr\"{u}fer associated with $\psi _{a_{k}}(\varphi )$%
, it can be shown (see \cite{Pearson, P1, KR})%
\begin{equation*}
\rho (\Sigma )=\lim_{k\rightarrow \infty }\rho _{k}(\Sigma
)=\lim_{k\rightarrow \infty }\frac{2}{\pi }\int_{\Sigma }\frac{\sin
^{2}\varphi }{R_{k}^{2}(\varphi )}d\varphi ~
\end{equation*}%
for any Borel set $\Sigma \subset \left( 0,\pi \right) $. Pearson's idea is
that sparse `barriers' lead to `independence' of certain (deterministic)
functions which behave as functions of an uniformly distributed random
variable. As a consequence, we have 
\begin{eqnarray}
\left( \frac{1}{R_{k}^{2}(\varphi )}\right) ^{1/k} &=&\prod_{m=1}^{k}\left( 
\frac{R_{m-1}^{2}(\varphi )}{R_{m}^{2}(\varphi )}\right) ^{1/k}  \notag \\
&=&\exp \left( \frac{1}{k}\sum_{m=1}^{k}\ln \frac{R_{m-1}^{2}(\varphi )}{%
R_{m}^{2}(\varphi )}\right)   \notag \\
&\equiv &\exp \left( \frac{1}{k}\sum_{m=1}^{k}\ln f\left( \varphi
,(a_{m}-a_{m-1})\varphi ,\theta _{m-1}(\varphi )\right) \right)   \notag \\
&\longrightarrow &\exp \left( \frac{1}{\pi }\int_{0}^{\pi }\ln f\left(
\varphi ,u,\theta \right) du\right) \equiv 1/r  \label{r}
\end{eqnarray}%
with probability one with respect to that uniform distribution, by the weak
law of large numbers. His method was modified in \cite{MarWre} by exploiting
the uniform distribution of a sequence $\left( \zeta _{m}(\varphi )\right)
_{m\geq 1}$, for almost every $\varphi $, defined by a linear interpolation
of Pr\"{u}fer angles:\footnote{%
Our definition of Pr\"{u}fer angles differs slightly from that of \cite%
{Pearson} and other authors. By $\theta _{j}$ we mean the Pr\"{u}fer angle
at the site $a_{j}$ immediately before the $j$--th barrier takes place.
Pearson's definition is at the point $b_{j}$ right after the barrier.}%
\begin{equation}
\theta _{m}(\varphi )=g(\theta _{m-1}(\varphi ))-(a_{m}-a_{m-1})\varphi
~,\qquad m\geq 2  \label{theta}
\end{equation}%
with $\theta _{1}=\theta _{0}-a_{1}\varphi ,$ for the monotone increasing
function 
\begin{equation*}
g(\theta )=\tan ^{-1}\left( (\tan \theta +\cot \varphi )/(1-\delta
)^{2}-\cot \varphi \right) 
\end{equation*}%
that maps the interval $(-\pi /2,\pi /2]$ into itself. The crucial
observation here is that $f(\varphi ,(a_{m}-a_{m-1})\varphi ,\theta
_{m-1}(\varphi ))$ in (\ref{r}) can be rewritten as $f(\varphi ,\theta
_{m}(\varphi ))$ for a different, although similar, function $f$. Equation (%
\ref{r}) thus gives an exact decay rate $1/r$ of $\psi _{n}(\varphi )$
without evoking `independence' of the Bernoulli shift sequence $%
u_{m}=(a_{m}-a_{m-1})\varphi ~\text{mod}~\pi $, which would require an
extremely sparse condition.

The Hausdorff dimension of the spectral measure $\rho $ can be determined
using an extension due to Jitormiskaya--Last \cite{JitLast} of the
Gilbert--Pearson theory of subordinance \cite{GP}, which relates the
spectral property of $\rho $ to the growth rate of solution $\psi _{n}^{\phi
}(\varphi )$ of the Schr\"{o}dinger equation $\Delta _{\delta ,\phi
}^{k}\psi _{n}=\lambda \psi _{n}$. Note that $\Delta _{\delta ,\phi
}^{k}=\Delta _{\delta ,0}^{k}+E_{0}\tan \phi $ and the phase boundary is
important since the exact Hausdorff dimension holds only for almost every $%
\phi $ w.r.t. the Lebesgue measure. It is worth mentioning that Zlato\v{s} 
\cite{Zla} has applied the Jitormiskaya--Last method to a sparse model very
similar to the one considered in \cite{MarWre} (whose `barriers' locate at
sites, not at edges). He has obtained the exact Hausdorff measure for a
sparse random model in which the distances from one to another `barrier' are
given by $a_{j}-a_{j-1}+\omega _{j}$ with $\left( \omega _{j}\right) _{j\geq
1}$ independent random variables uniformly distributed in the interval $%
\left[ -j,-j+1,\ldots ,j\right] $. The improvement of Pearson's method (\ref%
{r}) given in \cite{MarWre} allows the Hausdorff dimension of the spectral
measure to be determined without adding a random variable to the sparse
condition. This our main result (Theorem \ref{thethe}).

The present paper is organized as follows. We present some preliminary facts
on the spectrum of $\mathcal{J}_{\delta }$ on Section \ref{SFN} and on
Section \ref{EHD} we establish the exact Hausdorff dimension of the spectral
matrix measure of $\Delta _{\delta ,\phi }$. Our main result, Theorem \ref%
{thethe}, is stated and proved in this section, after we have extended to
the block--Jacobi matrix $\mathcal{J}_{\delta ,\phi }$ several preliminary
results of \cite{JitLast, Zla}.

\section{The Spectrum of $\mathcal{J}_{\protect\delta }$ and Notation\label%
{SFN}}

\setcounter{equation}{0} \setcounter{theorem}{0}

In order to introduce the spectral measure of block--Jacobi matrices
considered and to fix notation we shall first consider the $0$--Dirichlet
Laplacean operator $\Delta _{0.0}$. For convenience, we always change the
order of the tensor product in (\ref{kronecker}): $I_{L}\otimes J_{\delta
}+A_{L}\otimes I=\Pi \left( J_{\delta }\otimes I_{L}+I\otimes A_{L}\right)
\Pi ^{-1}$ by an appropriate permutation matrix $\Pi $ and we call it by $%
\mathcal{J}_{\delta }$ as well. \noindent The Kronecker sum 
\begin{equation}
\mathcal{J}_{0}=I_{L}\otimes J_{0}+A_{L}\otimes I\;,  \label{Pi}
\end{equation}%
with $J_{0}$ the free Jacobi matrix 
\begin{equation}
J_{0}=\left( 
\begin{array}{ccccc}
0 & 1 & 0 & 0 & \cdots \\ 
1 & 0 & 1 & 0 & \cdots \\ 
0 & 1 & 0 & 1 & \cdots \\ 
0 & 0 & 1 & 0 & \cdots \\ 
\vdots & \vdots & \vdots & \vdots & \ddots%
\end{array}%
\right) \;,  \label{JAZE}
\end{equation}%
\noindent is thus unitarily similar to (\ref{kronecker}) for (\ref{Pe}) with 
$\delta =0$; consequently, its spectrum remains unchanged.

The structure of (\ref{Pi}) permits to give a simple answer to the spectrum
of this operator. It is well known (see e.g. \cite{LANCASTER}) that if $%
\{\eta _{k}\}_{k=1}^{n}$ and $\{\lambda _{j}\}_{j=1}^{m}$ are the
eigenvalues of the matrices $A$ and $B$, respectively, then $\{\eta
_{k}+\lambda _{j}\}_{k,j=1}^{n,m}$ are the eigenvalues of the Kronecker sum $%
I_{m}\otimes A+B\otimes I_{n}$. Since the interval $[-2,2]$ is the essential
spectrum of $J_{0}$, the essential spectrum of $\mathcal{J}_{0}$ is given by 
\begin{equation}
\sigma _{\mathrm{ess}}(\mathcal{J}_{0})=\bigcup_{j=0}^{L-1}(\eta _{j}+\sigma
_{\mathrm{ess}}(J_{0}))\equiv \bigcup_{j=0}^{L-1}I_{j}\;,  \label{esesj}
\end{equation}%
\noindent with $\{\eta _{j}\}_{j=0}^{L-1}$, $\eta _{j}=2\cos \left( 2\pi
j/L\right) $, the eigenvalues of $A_{L}$. Thus, 
\begin{equation}
\sigma _{\mathrm{ess}}(\mathcal{J}_{0})=\left\{ 
\begin{array}{lll}
\left[ -2+2\cos \left( \pi (L-1)/L\right) ,4\right] & \mathrm{if} & L~%
\mathrm{is~odd} \\ 
\left[ -4,4\right] & \mathrm{if} & L~\mathrm{is~even}%
\end{array}%
\right.  \label{SUPO}
\end{equation}%
holds for $L\geq 2$.


It is also well known that the essential spectrum of the free Jacobi matrix $%
J_{0}$, defined by (\ref{JAZE}), is purely absolutely continuous. As $%
\mathcal{J}_{0}$ is in some sense a free matrix, we have

\begin{proposition}
\label{TEO3} The essential spectrum of $\mathcal{J}_{0}$, given by (\ref%
{SUPO}), is purely absolutely continuous.
\end{proposition}

\noindent \textit{Proof.} Let 
\begin{equation*}
M(z)=\int \frac{d\rho (x)}{x-z}
\end{equation*}%
be the $L\times L$ $M$--matrix defined by the Borel transform of the
spectral matrix $\rho $. By the spectral theorem, the $M$--matrix of $%
\mathcal{J}$ is related to the resolvent matrix $(\mathcal{J}-zI_{L}\otimes
I)^{-1}$ as follows. If $\mathcal{J}$ is the matrix representation of a
self--adjoint operator $H$ in the separable space $\mathcal{H}$ with an
orthonormal basis $\left\{ \varphi _{(k,m)}\right\} _{(k,m)\in \Lambda }$,
we have%
\begin{equation*}
(\mathcal{J}-zI_{L}\otimes I)_{(m,0)(m^{\prime },0)}^{-1}=(\varphi
_{(0,m)},\left( H-zI\right) ^{-1}\varphi _{(0,m^{\prime })})=\int \frac{%
d\rho _{mm^{\prime }}(x)}{x-z}=M_{mm^{\prime }}(x)
\end{equation*}

By the fact that $\Delta _{0,0}$ has periodic condition on the vertical
direction, $A_{L}$ is cyclic and the resolvent can be block--diagonalized by
the Fourier matrix: 
\begin{equation}
(F_{L}^{-1}\otimes I)(\mathcal{J}-zI_{L}\otimes I)^{-1}(F_{L}\otimes I)=%
\text{diag}\left\{ (J_{0}-z_{j}I)^{-1}\right\} _{j=0}^{L-1}\;,  \label{MaM}
\end{equation}%
\noindent with $z_{j}=z-2\cos (2\pi j/L)$, $j=0,\ldots ,L-1$, and $%
F_{L}:=[v_{1}v_{2}\cdots v_{L}]$ the matrix built up with the eigenvectors $%
v_{k}=(1,\xi ^{k},\ldots ,\xi ^{(L-1)k})/\sqrt{L}$, $\xi =\exp \{2\pi i/L\}$
of the shift matrix $S:\left( x_{0},\ldots ,x_{L-1}\right) \longrightarrow
\left( x_{1},\ldots ,x_{L-1},x_{0}\right) $ on its columns.

The $M$--matrix can thus be written as 
\begin{equation}
M(z)=F_{L}\text{diag}\left\{ (J_{0}-z_{j}I)_{00}^{-1}\right\}
_{j=1}^{L}F_{L}^{-1}\;,  \label{EME}
\end{equation}%
\noindent where the $00$--element $(J_{0}-z_{j}I)_{00}^{-1}$ of the
resolvent matrix $\left( J_{0}-z_{j}I\right) ^{-1}$ is the Weyl--Titchmarsh $%
m$--function of the free Jacobi matrix $J_{0}$ evaluated at $z_{j}$. It is a
simple exercise to calculate the $m$--function for the one-dimensional free
problem. If $u_{1}$, $u_{2}$ are the linear independent solutions of the Schr%
\"{o}dinger equation $J_{0}u=zu$ satisfying, respectively, Dirichlet (\ref{D}%
) and Neumann (\ref{N}) boundary conditions at $n=-1$, $m(z)$ is uniquely
defined by imposing that $u=u_{2}-m(z)u_{1}$ is $l^{2}(\mathbb{Z}_{+},%
\mathbb{C})$. Explicitly 
\begin{equation}
m(z_{j})=-\frac{z_{j}}{2}+\sqrt{\frac{z_{j}^{2}}{4}-1}\;,\qquad \qquad
j=0,1,\ldots ,L-1\;.  \label{eme}
\end{equation}%
\noindent for $\Re z_{j}>0$.

Now, let $m(z)=\displaystyle\int d\mu (x)/(z-x)$, 
\begin{equation*}
\Im m(\zeta )=\limsup_{\xi \downarrow 0}\Im m(z)\;,
\end{equation*}%
\noindent $z=\zeta +i\xi $, and let $L(\rho )$ be the set of all $\zeta \in 
\mathbb{R}$ for which this limit exists. It is known (see Appendix B from 
\cite{LASEIEU}) that the minimal (or essential) supports $\mathcal{M}$, $%
\mathcal{M}_{\mathrm{ac}}$ and $\mathcal{M}_{\mathrm{s}}$ of $\mu $, the
absolutely continuous part $\mu _{\mathrm{ac}}$ and the singular part $\mu _{%
\mathrm{s}}$ of $\mu $, with respect to the Lebesgue measure in $\mathbb{R}$%
, are, respectively, given by $\zeta \in L(\rho )$ such that $0<\Im
\,m(\zeta )\leq \infty $, $0<\Im \,m(\zeta )<\infty $ and $\Im \,m(\zeta
)=\infty $. These criteria can be obtained using de la Vall\'{e}e-Poussin's
decomposition theorem \cite{Saks}, the Lebesgue-Radon-Nikodym theorem and
the following Lemma (see e.g. \cite{GP}):

\begin{lemma}
\label{RadNik} If $(d\mu/d\nu)(\zeta)$ (the Radon-Nikodym derivative) exists
finitely or infinitely, them $\Im m(\zeta)$ also exists and $%
(d\mu/d\nu)(\zeta)=(1/\pi)\Im m(\zeta)$ ($\nu$ is some Lebesgue measure on $%
\mathbb{R}$).
\end{lemma}

Returning to the $M$--matrix (\ref{EME}), its diagonal elements are given by 
\begin{eqnarray}
M_{mm}(z) &=&\sum_{j,k=0}^{L-1}(F_{L})_{mj}\left( \text{diag}\left\{ \left(
J_{0}-z_{l}I\right) _{00}^{-1}\right\} _{l=1}^{L}\right) _{jk}\left(
F_{L}^{-1}\right) _{km}  \notag \\
&=&\sum_{j=1}^{L}m(z_{j})\left\vert (F_{L})_{mj}\right\vert ^{2}=\frac{1}{L}%
\sum_{j=1}^{L}m(z_{j})\;,  \label{FUe}
\end{eqnarray}%
\noindent with $m(z)$ given by (\ref{eme}). This equation, together with 
\begin{equation*}
\lim_{\xi \downarrow 0}\Im m(\zeta _{j}+i\xi )=\left\{ 
\begin{array}{c@{\quad \mathrm{if} \quad}l}
0 & \left\vert \zeta _{j}\right\vert \geq 2\, \\ 
\sqrt{1-\zeta _{j}^{2}/4} & \left\vert \zeta _{j}\right\vert <2\;%
\end{array}%
\right. ~,
\end{equation*}%
$\zeta _{j}=\zeta -2\cos (2\pi j/L)$, and Lemma \ref{RadNik}, leads to 
\begin{equation}
\lim_{\xi \downarrow 0}\frac{d\rho _{mm}}{d\zeta }(\zeta +i\xi )=\frac{1}{%
\pi L}\sum_{j=1}^{L}\lim_{\xi \downarrow 0}\Im m(\zeta _{j}+i\xi )\;,
\label{DENES}
\end{equation}%
\noindent which is strictly positive for almost every $\zeta $ with respect
to the Lebesgue measure on the essential support (\ref{SUPO}) of $\mathcal{J}%
_{0}$ and zero on its complement. The proof of Proposition \ref{TEO3} is
thus concluded evoking the above criteria.

\hfill $\Box $

A natural question to ask is whether the essential spectrum of the matrix $%
\mathcal{J}_{\delta }$ is, regardless of $\delta \in (0,1)$, the same of $%
\mathcal{J}_{0}$. This question is settled by the following

\begin{theorem}
\label{tee} Let $\mathcal{J}_{\delta }$ be the block--Jacobi matrix defined
by (\ref{kronecker}) with $P$ given by (\ref{Pe}) and (\ref{cesp}). The
essential spectrum of $\mathcal{J}_{\delta }$ is the set (\ref{esesj}) and,
consequently, 
\begin{equation*}
\sigma _{\mathrm{ess}}(\mathcal{J}_{\delta })=\sigma _{\mathrm{ess}}(%
\mathcal{J}_{0})
\end{equation*}%
holds for any $\delta \in $ $(0,1)$.
\end{theorem}

\begin{remark}
{Theorem \ref{tee} is an extension of Theorem 2.1 from \cite{MarWre}. We
follow its proof step by step.}
\end{remark}

\begin{remark}
The operator $\Delta _{\delta ,\phi }$ with $\phi $--phase boundary
condition at $k=-1$ (\ref{boucon}) may also be written in the block--Jacobi
matrix form (\ref{kronecker}) (see equation (\ref{JJ+})). Clearly $%
E_{0}\otimes \tan \phi I_{L}$ is a rank--$L$ perturbation of $\mathcal{J}%
_{\delta }$ and $\sigma _{\mathrm{ess}}(\mathcal{J}_{\delta ,\phi })=\sigma
_{\mathrm{ess}}(\mathcal{J}_{\delta })$, by Weyl's invariance principle (see
e.g. \cite{RESI}). Thus, it is sufficient to deal with $\mathcal{J}_{\delta
} $ to determine the essential spectrum of $\mathcal{J}_{\delta ,\phi }$.{\ 
\newline
}
\end{remark}

\noindent \textit{Proof.} Firstly, let us show that $\sigma _{\mathrm{ess}}(%
\mathcal{J}_{\delta })\subseteq \sigma _{\mathrm{ess}}(\mathcal{J}_{0})$.
Define for $u=(u(k,m))_{(k,m)\in \Lambda }\in l_{2}(\Lambda )$ the $2L$%
--dimensional column vectors 
\begin{equation*}
\mathbf{u}_{k}=(u(k,0),u(k+1,0),\ldots ,u(k,L-1),u(k+1,L-1))
\end{equation*}%
\noindent and the $2L\times 2L$ matrices 
\begin{equation}
h_{k}^{L}=p_{k}I_{L}\otimes A_{2}+\frac{1}{2}A_{L}\otimes I_{2}\;.
\label{AGA}
\end{equation}%
Then, the quadratic form associated with $\mathcal{J}_{\delta }$ can be
written as 
\begin{equation}
\left( u,\mathcal{J}_{\delta }u\right) =\sum_{k=1}^{\infty }\mathbf{u}%
_{k}\cdot h_{k}^{L}\mathbf{u}_{k+1}+\sum_{k=1}^{L}u(0,k)u(0,k+1)\;.
\label{AGGA}
\end{equation}%
The factor 1/2 present in (\ref{AGA}) avoids double counting of terms in (%
\ref{AGGA}); the second sum present in (\ref{AGGA}) corrects the counting of
the interacting terms between the elements of the first column.

We follow the strategy used in Proposition \ref{TEO3} to calculate the
eigenvalues of $h_{n}^{L}$. The characteristic polynomial of $h_{k}^{L}$
reads 
\begin{eqnarray*}
\det \left[ h_{k}^{L}-\lambda I_{L}\otimes I_{2}\right] &=&\det \left[
\left( F_{L}^{-1}\otimes I_{2}\right) \left( h_{k}^{L}-\lambda I_{L}\otimes
I_{2}\right) \left( F_{L}\otimes I_{2}\right) \right] \\
&=&\det \left[ \mathrm{diag}\left( S_{m}-\lambda I_{2}\right) _{m=0}^{L-1}%
\right] \\
&=&\prod_{m=0}^{L-1}\det \left[ S_{m}-\lambda I_{2}\right] \;,
\end{eqnarray*}%
\begin{equation*}
S_{m}=\left( 
\begin{array}{cc}
\cos \left( 2\pi m/L\right) & p_{k}\vspace{2mm} \\ 
p_{k} & \cos \left( 2\pi m/L\right)%
\end{array}%
\right) \;.
\end{equation*}%
So, the eigenvalues of $h_{k}^{L}$ are $\lambda _{k,m}^{\pm }=\pm p_{k}+\cos
(2\pi m/L)$, $m=0,\ldots ,L-1$. Inserting the spectral decomposition of $%
h_{k}^{L}$ 
\begin{equation*}
h_{k}^{L}=\sum_{m=0}^{L-1}\left( \lambda _{k,m}^{+}P_{k,m}^{+}+\lambda
_{k,m}^{-}P_{k,m}^{-}\right)
\end{equation*}%
\noindent into (\ref{AGGA}), where $P_{k,m}^{\pm }$ are the projectors in
the direction of the eigenvectors associated with $\lambda _{k,m}^{\pm }$,
we have 
\begin{equation*}
2\lambda ^{-}\leq \frac{\left( u,\mathcal{J}_{\delta }u\right) }{(u,u)}\leq
2\lambda ^{+}\;,
\end{equation*}%
\noindent with 
\begin{equation*}
\lambda ^{+}=\sup_{k,m}\lambda _{k,m}^{+}=2
\end{equation*}%
\noindent and 
\begin{equation*}
\lambda ^{-}=\inf_{k,m}\lambda _{k,m}^{-}=\left\{ 
\begin{array}{lll}
-2 & \mathrm{if} & L~\mathrm{is~even} \\ 
-1+\cos \left( \pi (L-1)/L\right) & \mathrm{if} & L~\mathrm{is~odd}%
\end{array}%
\right.
\end{equation*}%
concluding, together with (\ref{SUPO}), that $\sigma _{\mathrm{ess}}\left( 
\mathcal{J}_{\delta }\right) \subseteq \sigma _{\mathrm{ess}}\left( \mathcal{%
J}_{0}\right) $.

To prove the inclusion $\sigma _{\mathrm{ess}}\left( \mathcal{J}_{\delta
}\right) \supseteq \sigma _{\mathrm{ess}}\left( \mathcal{J}\right) $, we use
the Weyl criterion (Theorem VII.12 of \cite{RESI}): if $B$ is a bounded
self-adjoint operator on a separable Hilbert space $\mathscr{H}$, $\lambda $
belongs to the spectrum $\sigma (B)$ of $B$ if and only if there exists a
sequence $\left( \psi _{n}\right) _{n\in \mathbb{N}}$ in $\mathscr{H}$, with 
$\left\Vert \psi _{n}\right\Vert =1$, such that 
\begin{equation*}
\lim_{n\rightarrow \infty }\left\Vert (B-\lambda )\psi _{n}\right\Vert =0\;.
\end{equation*}%
Let $\lambda _{m}(\varphi )=2\left( \cos \varphi -\cos (2\pi m/L)\right) $, $%
\varphi \in \lbrack 0,\pi ]$, $m=0,\ldots ,L-1$, and define 
\begin{equation*}
\psi _{n,m}=\psi _{n}\otimes v_{m}\;,
\end{equation*}%
\noindent with $\psi _{n}=\left( 1/\sqrt{n}\right) (e^{i\varphi _{j}},\ldots
,e^{in\varphi _{j}},0,\ldots )$ and $v_{m}$ the $m$--th eigenvector of the
shift operator $S$ (see equation (\ref{MaM})). Clearly $\psi _{n,m}\in
l_{2}(\Lambda )$ and $\left\{ \lambda _{m}(\varphi )\right\} $ is in
one--to--one correspondence with (\ref{esesj}). We claim that, for each $%
m=0,\ldots ,L-1$, 
\begin{equation}
\left\Vert \mathcal{J}_{\delta }\psi _{n,m}-\;\lambda _{m}\psi
_{n,m}\right\Vert \leq c\frac{\ln n}{\sqrt{n}}  \label{NORMA}
\end{equation}%
\noindent holds with $c=c(\beta )$ independent of $n$. To prove (\ref{NORMA}%
), we just have to note that $\mathcal{J}_{\delta }\psi _{n,m}-\;\lambda
_{m}\psi _{n,m}$ consists of the action on $\psi _{n,m}$ of a sum of local
matrices, bounded in norm by one; $2$ of them involve the extreme points $%
e^{i\varphi }$ and $e^{in\varphi }$, and there are $\mathit{O}(\ln n)$
nondiagonal matrices. The $\mathit{O}(\ln n)$ is due to the fact that the
sequence $(a_{j})_{j\geq 1}$ satisfies the sparseness condition (\ref{cesp}%
), with at most $r$ points $a_{j}$ within $[1,n]$; $r$ is such that 
\begin{equation*}
r\leq \frac{\ln n}{\ln \beta }\;.
\end{equation*}%
Note that $A_{L}v_{m}=2\cos (2\pi m/L)v_{m}$ and this part of the tensor
product in $\mathcal{J}_{\delta }$ has no effect to the limit process. This
proves the inclusion $\sigma _{\mathrm{ess}}\left( \mathcal{J}_{\delta
}\right) \supseteq \sigma _{\mathrm{ess}}\left( \mathcal{J}\right) $ and
completes the proof of Theorem \ref{tee}.

\hfill $\Box $

\section{Exact Hausdorff dimension}

\label{EHD}

\setcounter{equation}{0} \setcounter{theorem}{0}

This section is devoted to the determination of the Hausdorff dimension of
the spectral measure of (\ref{kronecker}).

\subsection{Basic Definitions and Subordinacy}

We start by some useful definitions. A more complete description is found in 
\cite{Last}.

Given a Borel set $S\subset \mathbb{R}$ and $\alpha \in \lbrack 0,1]$, we
define the number 
\begin{equation}
Q_{\alpha ,\delta }(S)=\inf \left\{ \sum_{\nu =1}^{\infty }|b_{\nu
}|^{\alpha }:|b_{\nu }|<\delta ;S\subset \bigcup_{\nu =1}^{\infty }b_{\nu
}\right\} \;,  \label{Qhaus}
\end{equation}%
\noindent the infimum taken over all $\delta $--covers by intervals of size
at most $\delta $. The limit $\delta \rightarrow 0$, 
\begin{equation}
h^{\alpha }(S)=\lim_{\delta \downarrow 0}Q_{\alpha ,\delta }(S)\;,
\label{Mhaus}
\end{equation}%
\noindent is called $\alpha $-\textit{dimensional Hausdorff measure}. This
measure can be viewed as a continuous interpolation of the counting measure
at $\alpha =0$ (which assigns to each set $S$ the number of points in it)
and the Lebesgue measure at $\alpha =1$. It is clear by the definitions (\ref%
{Qhaus}) and (\ref{Mhaus}) that $h^{\alpha }(S)$ is an outer measure on $%
\mathbb{R}$, and its restriction to Borel sets is a Borel measure (see e.g. 
\cite{Falconer}). For $\beta <\alpha <\gamma $, 
\begin{equation*}
\delta ^{\alpha -\gamma }Q_{\gamma ,\delta }(S)\leq Q_{\alpha ,\delta
}(S)\leq \delta ^{\alpha -\beta }Q_{\beta ,\delta }(S)\;,
\end{equation*}%
\noindent holds for any $\delta >0$ and $S\subset \mathbb{R}$. So, if $%
h^{\alpha }(S)<\infty $, then $h^{\gamma }(S)=0$ for $\gamma >\alpha $; if $%
h^{\alpha }(S)>0$, then $h^{\beta }(S)=\infty $ for $\beta <\alpha $. Thus,
for every Borel set $S$, there is an unique $\alpha _{S}$ such that $%
h^{\alpha }(S)=0$ if $\alpha >\alpha _{S}$ and $h^{\alpha }(S)=\infty $ if $%
\alpha _{S}<\alpha $. The number $\alpha _{S}$ is called the \textit{%
Hausdorff dimension} of the set $S$.

Another useful concept is the exact dimension of a measure, due to
Rodgers-Taylor \cite{RodTay}:

\begin{definition}
\label{Exactdim} A measure $\mu $ defined on $\mathbb{R}$ is said to be of
exact dimension $\alpha $, $\alpha \in \lbrack 0,1]$, if and only if two
requirements hold: (1) for every $\beta \in \lbrack 0,1]$ with $\beta
<\alpha $ and $S$ a set of dimension $\beta $, $\mu (S)=0$ (which means that 
$\mu (S)$ gives zero weight to any set $S$ with $h^{\alpha }(S)=0$); (2)
there is a set $S_{0}$ of dimension $\alpha $ which supports $\mu $ in the
sense that $\mu (\mathbb{R}\backslash S_{0})=0$.
\end{definition}

Given a positive, finite measure $\mu$ and $\alpha \in [0,1]$, we define the
Hausdorff upper derivative by the limit 
\begin{eqnarray}
D_{\mu}^{\alpha}(x) \equiv \limsup_{\epsilon \downarrow 0} \frac{%
\mu((x-\epsilon, x+\epsilon))}{(2\epsilon)^{\alpha}}\;.  \label{UpDer}
\end{eqnarray}

Definition (\ref{UpDer}) is the generalization of the Radon-Nikodym
derivative for Hausdorff measures. Note that the limit $\epsilon \downarrow
0 $ does not need to be defined. Clearly, if $D_{\mu }^{\alpha
}(x_{0})<\infty $ for some $x_{0}$ then, for all $\beta <\alpha $, 
\begin{equation*}
D_{\mu }^{\beta }(x_{0})=\limsup_{\epsilon \downarrow 0}(2\epsilon )^{\alpha
-\beta }\frac{\mu ((x-\epsilon ,x+\epsilon ))}{(2\epsilon )^{\alpha }}%
=\limsup_{\epsilon \downarrow 0}(2\epsilon )^{\alpha -\beta }D_{\mu
}^{\alpha }(x)=0\;.
\end{equation*}

In a similar fashion, if $D_{\mu}^{\alpha}(x_0) > 0$ for some $x_0$, then $%
D_{\mu}^{\beta}(x_0) = \infty$ for all $\beta > \alpha$. Thus, we can define
for each $x_0$ the local Hausdorff dimension $\alpha(x_0)$, given by 
\begin{eqnarray}
\alpha_\mu(x_0) \equiv \liminf_{\epsilon \downarrow 0} \frac{%
\ln\mu((x-\epsilon, x+\epsilon))}{\ln(2\epsilon)}\;.  \label{LHD}
\end{eqnarray}

Finally, we introduce the notion of continuity and singularity of a measure
with respect to the Hausdorff measure. Given $\alpha \in \lbrack 0,1]$, a
measure $\mu $ is called $\alpha $--continuous if $\mu (S)=0$ for every set $%
S$ with $h^{\alpha }(S)=0$; it is called $\alpha $--singular if it is
supported on some set $S$ with $h^{\alpha }(S)=0$. We can reformulate
Definition \ref{Exactdim} in this context: a measure $\mu $ is said to have
exact dimension $\alpha $ if, for every $\epsilon >0$, it is simultaneously $%
(\alpha -\epsilon )$--continuous and $(\alpha +\epsilon )$--singular.

The following remarkable result is due to Rodgers-Taylor \cite{RodTay} and
was extracted from Del Rio-Jitomirskaya-Last-Simon \cite{DRJLS}:

\begin{theorem}[Rodgers-Taylor]
\label{RTdec} Let $\mu $ be any measure and $\alpha \in \lbrack 0,1]$. Let 
\begin{equation*}
T_{\infty }=\{x:D_{\mu }^{\alpha }(x)=\infty \}
\end{equation*}%
\noindent and let $\chi _{\alpha }$ denote its characteristic function. Let $%
d\mu _{\alpha s}=\chi _{\alpha }d\mu $ and $d\mu _{\alpha c}=(1-\chi
_{\alpha })d\mu $. Then $d\mu _{\alpha s}$ and $d\mu _{\alpha c}$ are,
respectively, singular and continuous with respect to $h^{\alpha }$.
\end{theorem}

\begin{remark}
{The restriction $\mu (T_{+}\cap \cdot )$ to the set $T_{+}=\{x:0<D_{\mu
}^{\alpha }(x)<\infty \}$ is absolutely continuous with respect to $%
h^{\alpha }$, in the sense that it is given by $f(x)dh^{\alpha }(x)$ for
some $f\in \mathcal{L}^{1}(\mathbb{R},dh^{\alpha })$.}
\end{remark}

\begin{remark}
Theorem \ref{RTdec} permits an extension of the standard Lebesgue
decomposition of a Borel measure into continuous and singular parts, with
respect to the Hausdorff measure. The decomposition into absolutely
continuous, singular-continuous and pure point parts can also be extended
(see \cite{Last} for a complete study). All these measure decompositions
lead to a corresponding spectral decomposition of the Hilbert space.
\end{remark}

Let $J$ be an essentially self--adjoint operator on $l^{2}(\mathbb{Z}_{+})$
given by a Jacobi matrix and let 
\begin{equation}
Ju=\lambda u\;,  \label{ScrEq}
\end{equation}%
be the corresponding Schr\"{o}dinger equation. Jitomirskaya-Last \cite%
{JitLast} extended, for Hausdorff measures, the Gilbert--Pearson theory of
subordinacy \cite{GP}, for Lebesgue measures, which relates the spectral
property of $\rho $ to the rate of growth of the solutions of the Schr\"{o}%
dinger equation. A solution $u$ of (\ref{ScrEq}) is said to be subordinate
if 
\begin{equation*}
\lim_{l\rightarrow \infty }\frac{\left\Vert u\right\Vert _{l}}{\left\Vert
v\right\Vert _{l}}=0
\end{equation*}%
\noindent holds for any linearly independent solution $v$ of (\ref{ScrEq}),
where $\left\Vert \cdot \right\Vert _{l}$ denotes the $l^{2}(\mathbb{Z}_{+})$%
--norm truncated at the length $l\in \mathbb{R}$, i.e., 
\begin{equation*}
\left\Vert u\right\Vert _{l}^{2}\equiv
\sum_{n=0}^{[l]}|u(n)|^{2}+(l-[l])|u([l]+1)|^{2}\;,
\end{equation*}%
\noindent $\lbrack l]$ the integer part of $l$.

We shall see that the theory in \cite{MarWre} permits to distinguish
different kinds of singular-continuous spectra, suitable for the study of
the spectral measure $\rho _{j}(\lambda )$ associated to each
one--dimensional component of $\mathcal{J}_{\delta }$, since their
singularity becomes more pronounced when $\lambda $ varies from the center
to the border of the spectrum (see Theorem 4.4 of \cite{MarWre}).

To extend the block-diagonalization ideas used in Section \ref{SFN} to study
the spectral measure of $\Delta _{\delta ,\phi }$, given by (\ref{Pip}), we
define operators 
\begin{equation*}
(H_{\delta ,\phi }^{j}\psi )(n)=p_{n}\psi (n+1)+p_{n-1}\psi (n-1)+V_{j}\psi
(n)\;,
\end{equation*}%
\noindent on $l^{2}\left( \mathbb{Z}_{+},\mathbb{C}\right) $ subjected to a $%
\phi $--boundary condition at $n=-1$: 
\begin{equation}
\psi (-1)\cos \phi -\psi (0)\sin \phi =0\;,  \label{boucon2}
\end{equation}%
\noindent for each $j\in \{0,\ldots ,L-1\}$. The \textquotedblleft
potential\textquotedblright\ $V_{j}=2\cos (2\pi j/L)$ arises from the
block-diagonalization of $\Delta _{\delta ,\phi }$ by the Fourier matrix $%
F_{L}\otimes I$. Note that each $H_{\delta ,\phi }^{j}$ is the projection of 
$\Delta _{\delta ,\phi }$ into its $j$--th one-dimensional constituent.

To each $H_{\delta ,\phi }^{j}$ there corresponds a Schr\"{o}dinger equation 
\begin{equation}
J_{\delta }u_{j}=\lambda _{j}u_{j}\;,  \label{ScrEqP}
\end{equation}%
\noindent with $J_{\delta }$ given by (\ref{MAJA}); we incorporate the
factor $V_{j}$ to the spectral parameter $\lambda $ and define $\lambda
_{j}=\lambda -2\cos (2\pi j/L)$.

Now, let $\lambda \in \mathbb{R}$ and $u_{1,j}$ be the solution of (\ref%
{ScrEqP}) which satisfies the Dirichlet boundary condition at $-1$, namely 
\begin{equation}
u_{1,j}(-1)=0,\qquad \qquad u_{1,j}(0)=1\;,  \label{D}
\end{equation}%
and let $u_{2,j}$ be the solution which satisfies the Neumann boundary
condition 
\begin{equation}
u_{2,j}(-1)=1,\qquad \qquad u_{2,j}(0)=0\;.  \label{N}
\end{equation}

Following Jitomirskaya-Last \cite{JitLast}, we define for any given $%
\epsilon >0$ and for each $j\in \{1,\ldots ,l\}$ a length $l_{j}(\epsilon
)\in (0,\infty )$ by the equality 
\begin{equation}
\left\Vert u_{1,j}\right\Vert _{l_{j}(\epsilon )}\left\Vert
u_{2,j}\right\Vert _{l_{j}(\epsilon )}=\frac{1}{2\epsilon }  \label{eps}
\end{equation}%
\noindent (see equation (1.12) from \cite{JitLast}).

Since at most one of the solutions $\left\{ u_{1,j},u_{2,j}\right\} $ of (%
\ref{ScrEqP}) is $l^{2}$ (thanks to the Wronskian constancy), the left-hand
side of (\ref{eps}) is a monotone increasing function of $l$ which vanishes
at $l=0$ and diverges as $l\rightarrow \infty $. On the other hand, the
right-hand side of (\ref{eps}) is a monotone decreasing function of $%
\epsilon $ which diverges as $\epsilon \rightarrow 0$. We conclude that the
function $l(\epsilon )$ is a well defined monotone decreasing and continuous
function of $\epsilon $ which diverges as $\epsilon \rightarrow 0$.

$l_{j}(\epsilon )$ being defined, we can apply Theorem 1.1 of \cite{JitLast}
for each Weyl-Titchmarsh $m_{j}$-function related to each pair of solutions $%
u_{1,j}$ and $u_{2,j}$, $j\in \{1,\ldots ,l\}$: for fixed $\epsilon >0$%
\begin{equation*}
\frac{5-\sqrt{24}}{m_{j}(\lambda +i\epsilon )}\leq \frac{\left\Vert
u_{1,j}\right\Vert _{l_{j}(\epsilon )}}{\left\Vert u_{2,j}\right\Vert
_{l_{j}(\epsilon )}}\leq \frac{5+\sqrt{24}}{m_{j}(\lambda +i\epsilon )}~.
\end{equation*}%
Note that Theorem 1.2 of \cite{JitLast} and its corollaries also holds: if $%
\mu _{j}$ denotes the spectral measure of $H_{\delta ,\phi }^{j}$, then,
with $b=\alpha /(2-\alpha )$,%
\begin{equation}
\limsup_{\varepsilon \rightarrow 0}\frac{\mu _{j}\left( (\lambda -\epsilon
,\lambda +\epsilon )\right) }{(2\varepsilon )^{\alpha }}=\infty   \label{sup}
\end{equation}%
if and only if%
\begin{equation}
\liminf_{l\rightarrow \infty }\frac{\left\Vert u_{1,j}\right\Vert _{l}}{%
\left\Vert u_{2,j}\right\Vert _{l}^{b}}=0~.  \label{inf}
\end{equation}

\subsection{Extension to Block--Jacobi Matrices}

We may ask whether these results can be extended to the diagonal elements $%
\rho _{mm}$ of the spectral matrix $\rho $ of $\Delta _{\delta ,\phi }$. The
generalization of Theorem 1.2 from \cite{JitLast} is as follows. Since all
diagonal elements of $M$ are equal, it is enough to consider $\rho _{00}$.

\begin{theorem}
\label{t1.2} Let $\Delta _{\delta ,\phi }$ be given by (\ref{Pip}), $\lambda
\in \mathbb{R}$ and $\alpha \in (0,1)$. Then 
\begin{equation}
D_{\rho _{00}}^{\alpha }(\lambda )=\limsup_{\epsilon \downarrow 0}\frac{\rho
_{00}(\lambda -\epsilon ,\lambda +\epsilon ))}{(2\epsilon )^{\alpha }}%
=\infty   \label{t1.2c}
\end{equation}%
\noindent if and only if 
\begin{equation}
\liminf_{l\rightarrow \infty }\frac{\left\Vert u_{1,j}\right\Vert _{l}}{%
\left\Vert u_{2,j}\right\Vert _{l}^{b}}=0  \label{t1.2b}
\end{equation}%
\noindent for at least one $j\in \mathcal{I}(\lambda )$, where 
\begin{equation}
\mathcal{I}(\lambda ):=\left\{ m\in \{0,\ldots ,L-1\}:I_{m}\ni \lambda
\right\} ~,  \label{Ilambda}
\end{equation}%
$I_{m}$ is defined in (\ref{esesj}) and $b=\dfrac{\alpha }{2-\alpha }$.
\end{theorem}

\noindent \textit{Proof.} Suppose that (\ref{t1.2c}) holds. Then, by (\ref%
{FUe}) and (\ref{DENES}), there exists at least one $j\in \mathcal{I}%
(\lambda )$ such that 
\begin{equation}
\limsup_{\epsilon \downarrow 0}\frac{\mu _{j}(\lambda -\epsilon ,\lambda
+\epsilon ))}{(2\epsilon )^{\alpha }}=\infty ~  \label{t1.2a}
\end{equation}%
and by Theorem 1.2 in \cite{JitLast} applied to the operator $H_{\delta
,\phi }^{j}$ (equations (\ref{sup}) and (\ref{inf})), this holds if and only
if (\ref{t1.2b}) holds.

Suppose now that (\ref{t1.2b}) holds for some $j\in \mathcal{I}(\lambda )$.
The same Theorem 1.2 of \cite{JitLast} leads to (\ref{t1.2a}). But we know
from (\ref{FUe}) and (\ref{DENES}) that this implies (\ref{t1.2c}),
concluding the proof of Lemma \ref{t1.2}.\newline

\hfill $\Box $

The resulting corollaries of Theorem 1.2 in \cite{JitLast} can be extended
on a similar fashion. Of particular interest are Corollaries 4.4 and 4.5 of 
\cite{JitLast}. The new version of the first is given by

\begin{corollary}
\label{C4.4} Suppose that for some $\alpha \in \lbrack 0,1)$ and every $%
\lambda $ in some Borel set $A$, every solution $v_{j}$ of (\ref{ScrEq})
obeys 
\begin{equation*}
\limsup_{l\rightarrow \infty }\frac{\left\Vert v_{j}\right\Vert _{l}^{2}}{%
l^{2-\alpha }}<\infty
\end{equation*}%
\noindent for all $j\in \mathcal{I}(\lambda )\neq \emptyset $. Then the
restriction $\rho _{00}(A\cap \cdot )$ is $\alpha $-continuous.
\end{corollary}

\noindent \textit{Proof.} The proof follows the same structure of the proof
of Corollary 4.4 in \cite{JitLast}. Let $\lambda \in A$. From the constancy
of the Wronskian, $\left\Vert u_{1,j}\right\Vert _{l}\left\Vert
u_{2,j}\right\Vert _{l}\geq l$ holds for every $j$, and since, by
hypothesis, $\left\Vert u_{2,j}\right\Vert _{l}^{2}<Cl^{2-\alpha }$ for some
constant $C$, it follows that $\left\Vert u_{1,j}\right\Vert
_{l}>C^{-1/2}l^{\alpha /2}$ for every $j\in \mathcal{I}(\lambda )$. Thus, we
have 
\begin{equation*}
\frac{\left\Vert u_{1,j}\right\Vert _{l}}{\left\Vert u_{2,j}\right\Vert
_{l}^{b}}>C^{-(1+b)/2}l^{\alpha /2-b(2-\alpha )/2}=C^{-(1+b)/2}>0\;,
\end{equation*}%
\noindent since $b=\alpha /(2-\alpha )$. It follows from Theorem \ref{t1.2}
that $\rho _{00}(A\cap \cdot )$ is $\alpha $-continuous. \newline

\hfill $\Box $

Corollary \ref{C4.4} can be rewritten in terms of the one--dimensional $%
2\times 2$ transfer matrices 
\begin{equation}
T_{j}(n;\lambda )=T_{j}(n,n-1;\lambda )T_{j}(n-1,n-2;\lambda )\cdots
T_{j}(0,-1;\lambda )~,  \label{TTT}
\end{equation}%
where 
\begin{equation}
T_{j}(n,n-1;\lambda )=\left( 
\begin{array}{cc}
\displaystyle\frac{\lambda _{j}}{p_{n}} & \displaystyle\frac{-p_{n-1}}{p_{n}}%
\vspace{2mm} \\ 
1 & 0%
\end{array}%
\right) \equiv T(n,n-1;\lambda _{j})  \label{mt}
\end{equation}%
\noindent is related to the equation (\ref{ScrEqP}) for every $j\in
\{0,\cdots ,L-1\}$. Note that $T(n,n-1;\lambda )$ is precisely the transfer
matrix considered in \cite{MarWre} (see equation (2.2) therein). Moreover,
for the sequence $\left( p_{n}\right) _{n\geq -1}$ of the form (\ref{Pe}),
only three different $2\times 2$ matrices appear in the r.h.s. of (\ref{TTT}%
):%
\begin{equation}
T_{-}=\left( 
\begin{array}{cc}
\dfrac{\lambda _{j}}{1-\delta } & \dfrac{-1}{1-\delta } \\ 
1 & 0%
\end{array}%
\right) \ ,\qquad T_{+}=\left( 
\begin{array}{cc}
\lambda _{j} & -1+\delta \\ 
1 & 0%
\end{array}%
\right) \ \mathrm{and}\qquad T_{0}=\left( 
\begin{array}{cc}
\lambda _{j} & -1 \\ 
1 & 0%
\end{array}%
\right)  \label{T-T+T0}
\end{equation}%
depending on whether the left, the right or none of the two entries $n$ and $%
n-1$ in (\ref{mt}) belong to $\mathcal{A}$, respectively. As 
\begin{equation*}
\left( 
\begin{array}{c}
u_{j}(n+1)\vspace{2mm} \\ 
u_{j}(n)%
\end{array}%
\right) =T_{j}(n;\lambda )\left( 
\begin{array}{c}
u_{j}(0)\vspace{2mm} \\ 
u_{j}(-1)%
\end{array}%
\right) \;,
\end{equation*}%
$T_{j}(n;\lambda )$ is also the fundamental matrix of (\ref{ScrEqP}) 
\begin{equation}
T_{j}(n;\lambda )=\left( 
\begin{array}{cc}
u_{1,j}(n+1) & u_{2,j}(n+1)\vspace{2mm} \\ 
u_{1,j}(n) & u_{2,j}(n)%
\end{array}%
\right) \;.  \label{TraMa}
\end{equation}

Marchetti et al. \cite{MarWre} have determined precisely the growth of the
norm of $T(n;\lambda )$ given by the product of (\ref{mt}) with $\lambda
_{j}=\lambda $ and $P$ given by (\ref{Pe}) and (\ref{cesp}). This together
with a result due to Zlato\v{s} \cite{Zla} permits the determination of the
Hausdorff dimension of $\rho _{00}$.

Given (\ref{TraMa}), we have

\begin{corollary}
\label{C4.4a} Suppose that for some $\alpha \in \lbrack 0,1)$ and every $%
\lambda $ in some Borel set $A$, 
\begin{equation}
\limsup_{l\rightarrow \infty }\frac{1}{l^{2-\alpha }}\sum_{n=0}^{l}\left%
\Vert T_{j}(n;\lambda )\right\Vert ^{2}<\infty \;,  \label{dnor}
\end{equation}%
\noindent for all $j\in \mathcal{I}(\lambda )$, with $\left\Vert \cdot
\right\Vert $ some matrix norm. Then the restriction $\rho _{00}(A\cap \cdot
)$ is $\alpha $--continuous.
\end{corollary}

\noindent Proof. Theorem 2.3 from \cite{KLS} states that there are two
positive constants $c_{1}$, $c_{2}$, such that 
\begin{equation}
c_{1}\max \left\{ |u_{1,j}(n+1)|^{2},|u_{2,j}(n+1)|^{2}\right\} \leq
\left\Vert T_{j}(n;\lambda )\right\Vert ^{2}\leq c_{2}\max \left\{
u_{1,j}(n+1)|^{2},|u_{2,j}(n+1)|^{2}\right\}   \label{lr}
\end{equation}%
This leads to 
\begin{equation}
\sum_{n=0}^{l}\left\Vert T_{j}(n;\lambda )\right\Vert ^{2}\geq c\max
\{\left\Vert u_{1,j}\right\Vert _{l+1}^{2},\left\Vert u_{2,j}\right\Vert
_{l+1}^{2}\}  \label{c4.4a}
\end{equation}%
for every $j\in \mathcal{I}(\lambda )$. Hypothesis (\ref{dnor}), together
with (\ref{c4.4a}), implies Corollary \ref{C4.4a}. \newline

\hfill $\Box $

It is interesting to note that the growth of the norm of the transfer matrix
gives exactly the growth of the increasing solution. This fact will be of
great importance later.

The new version of Corollary 4.5 is

\begin{corollary}
\label{C4.5} Suppose that for at least one $j\in \mathcal{I}(\lambda )\neq
\emptyset $ 
\begin{equation}
\liminf_{l\rightarrow \infty }\frac{\left\Vert u_{1,j}\right\Vert _{l}^{2}}{%
l^{\alpha }}=0  \label{c4.5a}
\end{equation}%
\noindent for every $\lambda $ in some Borel set $A$. Then the restriction $%
\rho _{00}(A\cap \cdot )$ is $\alpha $--singular.
\end{corollary}

\noindent \textit{Proof.} Let $\lambda \in A$ and $b=\alpha /(2-\alpha )$.
By hypothesis, there is at least one $j\in \mathcal{I}(\lambda )$ that
satisfies (\ref{c4.5a}). Again, by the constancy of the Wronskian, $%
\left\Vert u_{1,j}\right\Vert _{l}\left\Vert u_{2,j}\right\Vert _{l}\geq l$,
and so $\left\Vert u_{2,j}\right\Vert _{l}^{b}\geq (l/\left\Vert
u_{1,j}\right\Vert _{l})^{b}$. This implies 
\begin{equation*}
\liminf_{l\rightarrow \infty }\frac{\left\Vert u_{1,j}\right\Vert _{l}}{%
\left\Vert u_{2,j}\right\Vert _{l}^{b}}\leq \liminf_{l\rightarrow \infty }%
\frac{\left\Vert u_{1,j}\right\Vert _{l}^{1+b}}{l^{b}}=\liminf_{l\rightarrow
\infty }\left( \frac{\left\Vert u_{1,j}\right\Vert _{l}^{2}}{l^{\alpha }}%
\right) ^{1/(2-\alpha )}=0\;.
\end{equation*}%
It follows from Theorem \ref{t1.2} that $\rho _{00}(A\cap \cdot )$ is $%
\alpha $-singular. \newline

\hfill $\Box $

\subsection{Main Result}

In order to state the result concerning the Hausdorff dimension of the
measure $\rho _{00}$, we need a result due to Zlato\v{s} \cite{Zla} on the
growth and decay of the solutions of (\ref{ScrEqP}) in the $\text{span}%
\left\{ u_{1,j},u_{2,j}\right\} $. We shall give an improved version of
Lemma 2.1 of \cite{Zla}.

\begin{proposition}
\label{l2.1} Let $\mathcal{A=}\left( a_{n}\right) _{n\geq 1}$ be given by (%
\ref{cesp}), $\lambda \in \mathbb{R}$ and let us assume that, for $j\in 
\mathcal{I}(\lambda )$, the sequence $\left( \theta _{n}^{j}\right) _{n\geq
0}$ of Pr\"{u}fer angles, defined by (\ref{theta}) with $\varphi $ replaced
by $\varphi _{j}$, is uniformly distributed $\text{mod}$ $\pi $ for every $%
\theta _{0}^{j}\in \lbrack 0,\pi ]$ and almost every $\varphi _{j}\in
\lbrack 0,\pi ]$ (w.r.t. Lebesgue measure) where $2\cos \varphi _{j}=\lambda
_{j}=\lambda -2\cos (2\pi j/L)$. Then, there is a generalized eigenfunction $%
u_{j}$ (i.e., $u_{j}$ satisfies (\ref{ScrEqP}) and the phase boundary
condition (\ref{boucon2})) for energy $\lambda $ such that 
\begin{equation}
C_{n}^{-1}r_{j}^{n/2}\leq \left\vert u_{j}(a_{n}+1)\right\vert \leq
C_{n}r_{j}^{n/2}\;,  \label{lru}
\end{equation}%
holds for a constants $r_{j}>1$ given by (with $p=1-\delta $) 
\begin{equation}
r_{j}=r(p,\lambda _{j})=1+\frac{(1-p^{2})^{2}}{p^{2}(4-\lambda _{j}^{2})}
\label{rj}
\end{equation}%
and $C_{n}^{1/n}\searrow 1$ as $n\rightarrow \infty $. In addition, there
exists a subordinate solution $v_{j}$ for energy $\lambda $ such that, for
all sufficiently large $n$, 
\begin{equation}
\left\vert v_{j}(a_{n}+1)\right\vert \leq \tilde{C}_{n}r_{j}^{-n/2}\;
\label{dj}
\end{equation}%
holds with $\tilde{C}_{n}^{1/n}\searrow 1$ as $n\rightarrow \infty $.
\end{proposition}

\noindent \textit{Proof.} We shall combine ideas of \cite{Zla} with Theorem
8.1 of \cite{LS} and estimates of \cite{MarWre}. Let us denote the spectral
norm of the transfer matrix $\left\Vert T_{j}(a_{n}+1;\lambda )\right\Vert $
by $t_{j,n}.$ Equation (\ref{lru}), together with (\ref{lr}), implies that $%
t_{j,n}$ satisfies the same upper and lower bounds. Under the hypotheses of
Proposition \ref{l2.1}, it follows from (3.8) and (4.19) of \cite{MarWre}
that 
\begin{equation}
C_{n}^{-1}r_{j}^{n/2}\leq t_{j,n}\leq C_{n}r_{j}^{n/2}~,  \label{rtr}
\end{equation}%
with $r_{j}$ given by (\ref{rj}).

By (\ref{mt}) and (\ref{Pe}), 
\begin{equation*}
\left\Vert T_{j}(k,k-1;\lambda )\right\Vert ^{2}\leq \left\Vert
T_{j}(k,k-1;\lambda )\right\Vert _{E}^{2}\leq 1+\frac{1+\lambda _{j}^{2}}{%
(1-\delta )^{2}}<\infty 
\end{equation*}%
if $\delta \in (0,1)$, where $\left\Vert \cdot \right\Vert _{E}$ is the
Euclidean matrix norm, for $k,k-1\in \mathcal{A}$; otherwise $%
T_{j}(k,k-1;\lambda )$ is similar to a clockwise rotation $R(\varphi _{j})$
by $\varphi _{j}=(1/2)\arccos \lambda _{j}$: $R(\varphi _{j})=UT_{0}U^{-1}$
(see (2.8) of \cite{MarWre}). We write 
\begin{equation*}
T_{j}(a_{n}+1;\lambda )=A_{n}\cdots A_{1}
\end{equation*}%
where, for each $m$ 
\begin{equation*}
A_{m}=T_{j}(a_{m}+1,a_{m};\lambda )\cdots T_{j}(a_{m-1}+2,a_{m-1}+1;\lambda
)=T_{-}T_{+}T_{0}^{\beta ^{m}-2}
\end{equation*}%
by (\ref{T-T+T0}). Denoting $s_{j,n}=$ $\left\Vert A_{n}\right\Vert $, we
thus have%
\begin{equation}
s_{j,n}\leq C\left( 1+\frac{1+\lambda _{j}^{2}}{(1-\delta )^{2}}\right)
\equiv B_{j}  \label{s}
\end{equation}%
$C=(1+\left\vert \cos \varphi _{j}\right\vert )/(1-\left\vert \cos \varphi
_{j}\right\vert )$, uniformly in $n$. As a consequence,%
\begin{equation}
\sum_{n=1}^{\infty }\frac{s_{j,n+1}^{2}}{t_{j,n}^{2}}<\infty   \label{sum}
\end{equation}%
verifies the assumption of Theorem 8.1 of \cite{LS} and provides the
existence of a subordinate solution $v_{j}$ for energy $\lambda $. The idea
of Zlato\v{s} is to use the proof of Last--Simon to establish the decay of
the subordinate solution. We shall reproduce the main steps, for convenience.

Since $T_{0}$, $T_{+-}:=T_{+}T_{-}$ given by (\ref{T-T+T0}) and,
consequently, $T_{j}(a_{n}+1;\lambda )$ and $T_{j}^{\ast }(a_{n}+1;\lambda )$
are $2\times 2$ unimodular real matrices, $T_{j}^{\ast }(a_{n}+1;\lambda
)T_{j}(a_{n}+1;\lambda )$ is a $2\times 2$ unimodular symmetric real matrix
whose eigenvalues are $t_{j,n}^{2}$ and $t_{j,n}^{-2}$, with corresponding
orthonormal eigenvectors $\mathbf{v}_{j,n}^{+}$ and $\mathbf{v}_{j,n}^{-}$: $%
\left( \mathbf{v}_{j,n}^{+},\mathbf{v}_{j,n}^{-}\right) =0$. We write $%
\mathbf{v}_{\alpha }=\left( 
\begin{array}{c}
\cos \alpha  \\ 
\sin \alpha 
\end{array}%
\right) $ and define $\alpha _{n}$ by 
\begin{equation}
\mathbf{v}_{\alpha _{n}}=\mathbf{v}_{j,n}^{-}~.  \label{vv}
\end{equation}%
Clearly, $\mathbf{v}_{j,n}^{+}=\mathbf{v}_{\alpha _{n}+\pi /2}$ and by the
spectral theorem, we have%
\begin{eqnarray}
\left\Vert T_{j}(a_{n}+1;\lambda )\mathbf{v}_{\alpha }\right\Vert ^{2}
&=&\left( \mathbf{v}_{\alpha },T_{j}^{\ast }(a_{n}+1;\lambda
)T_{j}(a_{n}+1;\lambda )\mathbf{v}_{\alpha }\right)   \notag \\
&=&t_{j,n}^{2}\left\vert (\mathbf{v}_{\alpha },\mathbf{v}_{+})\right\vert
^{2}+t_{j,n}^{-2}\left\vert (\mathbf{v}_{\alpha },\mathbf{v}_{-})\right\vert
^{2}  \notag \\
&=&t_{j,n}^{2}\sin ^{2}\left( \alpha -\alpha _{n}\right) +t_{j,n}^{-2}\cos
^{2}\left( \alpha -\alpha _{n}\right) ~.  \label{Ttt}
\end{eqnarray}

By the properties of a matrix norm together with (\ref{Ttt}) for $n+1$ and
definition (\ref{vv}), it can be shown (see proof of Theorem 8.1 of \cite{LS}%
)%
\begin{equation*}
\left\vert \alpha _{n}-\alpha _{n+1}\right\vert \leq \frac{\pi }{2}\frac{%
s_{j,n+1}^{2}}{t_{j,n}^{2}}~.
\end{equation*}

Condition (\ref{sum}) implies that the sequence $\left( \alpha _{n}\right)
_{n\geq 1}$ has a limit $\alpha ^{\ast }=\lim_{n\rightarrow \infty }\alpha
_{n}$. Hence, equation (\ref{Ttt}) and the telescope estimate%
\begin{equation*}
\left\vert \alpha _{n}-\alpha ^{\ast }\right\vert \leq \sum_{m=n}^{\infty
}\left\vert \alpha _{m}-\alpha _{m+1}\right\vert \leq \frac{\pi }{2}%
\sum_{m=n}^{\infty }\frac{s_{j,m+1}^{2}}{t_{j,m}^{2}}
\end{equation*}%
yields%
\begin{eqnarray*}
\left\Vert T(a_{n}+1;\lambda _{j})\mathbf{v}_{\alpha ^{\ast }}\right\Vert
^{2} &\leq &t_{j,n}^{2}\left( \alpha ^{\ast }-\alpha _{n}\right)
^{2}+t_{j,n}^{-2} \\
&\leq &\frac{\pi }{2}B_{j}t_{j,n}^{2}\left( \sum_{m=n}^{\infty }\frac{1}{%
t_{j,m}^{2}}\right) ^{2}+t_{j,n}^{-2}
\end{eqnarray*}%
which, together with (\ref{rtr}), gives (\ref{dj}) concluding the proof of
Proposition \ref{l2.1}. Note that, by definition of transfer matrix, $%
v_{j}(a_{n}+1)=\left( T_{j}(a_{n}+1;\lambda )\mathbf{v}_{\alpha ^{\ast
}}\right) _{2}$ is a subordinate solution evaluated at $a_{j}+1$ since $%
u_{j}(a_{n}+1)\equiv \left( T_{j}(a_{n}+1;\lambda )\mathbf{v}_{\alpha ^{\ast
}+\pi /2}\right) _{2}$ satisfies 
\begin{equation*}
\lim_{n\rightarrow \infty }\frac{\left\vert v_{j}(a_{n}+1)\right\vert }{%
\left\vert u_{j}(a_{n}+1)\right\vert }=0
\end{equation*}%
in view of $\left\Vert T(a_{n}+1;\lambda _{j})\mathbf{v}_{\alpha ^{\ast
}+\pi /2}\right\Vert \geq t_{j,n}^{2}/2$ for sufficiently large $n$.

\hfill $\Box $

\begin{remark}
Equation (\ref{rtr}), where $t_{j,n}\equiv \left\Vert T_{j}(a_{n}+1;\lambda
)\right\Vert $, holds for every $j\in \mathcal{I}(\lambda )$ for $\lambda
_{j}=\lambda -2\cos (2\pi j/L)\in (-2,2)\setminus A_{\theta _{0}^{j}}$, $%
A_{\theta _{0}^{j}}$ a set of zero Lebesgue measure possibly depending on
the initial Pr\"{u}fer angle $\theta _{0}^{j}$, which depends on $\phi $%
--condition and $\varphi _{j}$ (see eq. (3.7) and Theorem 4.4 of \cite%
{MarWre}).
\end{remark}

We are now ready to present our main result.

\begin{theorem}
\label{thethe} Let $\Delta _{\delta ,\phi }$ be given by (\ref{Pip}) with $%
\delta \in \left( 0,1\right) $ and $\phi $--boundary condition (\ref{boucon}%
). Let $\rho $ be its spectral matrix measure. For any closed interval of
energies $I\subset \bigcup_{j}I_{j}$, where 
\begin{equation}
I_{j}=\left( -2+2\cos \left( \frac{2\pi j}{L}\right) ,2+2\cos \left( \frac{%
2\pi j}{L}\right) \right) ~,  \label{INTI}
\end{equation}%
\noindent and for almost every boundary condition $\phi $, the element $\rho
_{00}$ of the spectral measure $\rho $ restricted to $I$ has, for every $%
\varepsilon >0$, the Hausdorff dimension%
\begin{equation}
\alpha _{\rho _{00}}(\lambda )\in \left( \alpha _{\rho _{j^{\ast }}}(\lambda
_{j^{\ast }})-\varepsilon ,\alpha _{\rho _{j^{\ast }}}(\lambda _{j^{\ast
}})+\varepsilon \right)   \label{epsilon}
\end{equation}%
where 
\begin{equation}
\alpha _{\rho _{j^{\ast }}}(\lambda _{j^{\ast }})=\min_{j\in \mathcal{I}%
(\lambda )}\alpha _{\rho _{j}}(\lambda _{j})=\min_{j\in \mathcal{I}(\lambda
)}\left( 1-\frac{\ln r_{j}}{\ln \beta }\right) \;,  \label{HDAL}
\end{equation}%
\noindent with $r_{j}=r(p,\lambda _{j})$ given by (\ref{rj}) ($p=1-\delta )$%
, if the sparseness parameter satisfies $\beta >\beta _{0}$ for some $\beta
_{0}=\beta _{0}(\delta ,\lambda _{j^{\ast }},\varepsilon )$ large enough.
\end{theorem}

\begin{remark}
{Theorem \ref{thethe} generalizes (from the one-dimensional case to the
finite strip problem) and improves (it establishes the Hausdorff dimension
of the spectral measure) Theorem 4.1 of Zlat\u{o}s \cite{Zla}. \newline
}
\end{remark}

\noindent \textit{Proof.} Let $I$ be given by (\ref{INTI}) and let us,
provisionally, assume that for $\lambda \in I$ the sequence $\left( \theta
_{n}^{j}\right) _{n\geq 0}$ of Pr\"{u}fer angles is uniformly distributed $%
\text{mod}$ $\pi $ for every $\theta _{0}^{j}\in \lbrack 0,\pi ]$ and almost
every $\varphi _{j}=\left( \cos ^{-1}\lambda _{j}\right) /2\in \lbrack 0,\pi
]$, for every $j\in \mathcal{I}(\lambda )$. It follows from (\ref{rtr}) and
Theorem 4.4 of \cite{MarWre} that, there is an $A_{\theta _{0}^{j}}$ with
zero Lebesgue measure such that for any $\lambda \in I\setminus A_{\theta
_{0}^{j}}$ and any $k\in \mathbb{Z}_{+}$ such that $a_{n}\leq k<a_{n+1}$, we
have 
\begin{equation*}
\left\Vert T_{j}(k;\lambda )\right\Vert \leq C_{n}r_{j}^{n/2}\leq
C_{n}^{\prime }a_{n}^{\gamma _{j}/2}\leq C_{n}^{\prime \prime }k^{\gamma
_{j}/2}\;,
\end{equation*}%
\noindent with $\gamma _{j}\equiv \ln r_{j}/\ln \beta $ and $%
\lim_{n\rightarrow \infty }\left( C_{n}^{\prime \prime }\right) ^{1/n}=1$,
by the sparseness condition (\ref{cesp}).

It follows from the constancy of $\left\Vert T_{j}(k;\lambda )\right\Vert $
on $[a_{n}+1,a_{n+1}]$ (see Section 4 of \cite{MarWre}), together with the
above equation, 
\begin{equation}
\sum_{k=0}^{l}\left\Vert T_{j}(k;\lambda )\right\Vert ^{2}\leq cl^{1+\gamma
_{j}}  \label{somatra}
\end{equation}%
\noindent\ holds for some $c>0$ and every $\lambda \in I\setminus A_{\theta
_{0}^{j}}$.

The application of Proposition \ref{l2.1} for these values of $\lambda $
guarantees the existence of a subordinate solution $v_{j}$ which satisfies 
\begin{equation*}
|v_{j}(a_{n}+1)|^{2}\leq C_{n}^{\prime \prime \prime }a_{n}^{-\gamma _{j}}
\end{equation*}%
\noindent for every $j\in \mathcal{I}(\lambda )$. Since every solution of (%
\ref{ScrEq}) has constant modulus on the interval $[a_{n}+1,a_{n+1}]$, we
have 
\begin{equation}
\left\Vert v_{j}\right\Vert _{l}^{2}\leq c^{\prime }l^{1-\gamma _{j}}\;,
\label{modusub}
\end{equation}%
\noindent for some $c^{\prime }>0$.

Since the measure $\rho _{00}$ restricted to $I$ is supported on the set of
those $\lambda $ for which each $u_{j}^{\mathrm{sub}}$ satisfies the
boundary condition $\phi $ (due to the fact that each constituent of $\rho
_{00}$ has no absolutely continuous part; see Theorem 1 of \cite{GP}), we
have $u_{1,j}=v_{j}$.

Thus, by (\ref{somatra}) and (\ref{modusub}) 
\begin{equation}
\limsup_{l\rightarrow \infty }\frac{1}{l^{2-\alpha }}\sum_{k=0}^{l}\left%
\Vert T_{j}(k;\lambda )\right\Vert ^{2}<\infty  \label{c4.4c}
\end{equation}%
\noindent and 
\begin{equation}
\liminf_{l\rightarrow \infty }\frac{\left\Vert u_{1,j}\right\Vert _{l}^{2}}{%
l^{\alpha ^{\prime }}}=0  \label{c4.5b}
\end{equation}%
\noindent hold for each $j\in \mathcal{I}(\lambda )$, provided $2-\alpha
\geq 1+\gamma _{j}$ and $\alpha ^{\prime }>1-\gamma _{j}$.

Corollary \ref{C4.4a} says that if (\ref{c4.4c}) is satisfied for all $j\in 
\mathcal{I}(\lambda )$, the restriction $\rho _{00}((I\setminus \cup
_{j}A_{\theta _{0}^{j}})\cap \cdot )$ is $\alpha $-continuous. Clearly, $%
\alpha =\min_{j}(1-\gamma _{j})$ satisfies the requirement: 
\begin{equation*}
\limsup_{l\rightarrow \infty }\frac{1}{l^{2-\alpha }}\sum_{n=1}^{l}\left%
\Vert T_{j}(n;\lambda )\right\Vert ^{2}\leq \limsup_{l\rightarrow \infty }%
\frac{1}{l^{1+\gamma _{j}}}\sum_{n=1}^{l}\left\Vert T_{j}(n;\lambda
)\right\Vert ^{2}<\infty \;,
\end{equation*}%
\noindent which implies that (\ref{c4.4c}) holds simultaneously for every $%
j\in \mathcal{I}(\lambda )$, provided $\lambda \in I\setminus \cup
_{j}A_{\theta _{j}}$. Thus $\rho _{00}((I\setminus \cup _{j}A_{\theta
_{0}^{j}})\cap \cdot )$ is at most $\alpha $-continuous.

We affirm that $\rho _{00}((I\setminus \cup _{j}A_{\theta _{0}^{j}})\cap
\cdot )$ is at least $\alpha $--singular with $\alpha =\min_{j}(1-\gamma
_{j})$. We have from Corollary \ref{C4.5} that the restriction above is $%
\eta $--singular for every $\eta >\alpha $ (since (\ref{c4.5a}) is satisfied
for at least one $j$). However, (\ref{c4.5b}) is satisfied for every $j$;
this proves our assertion.

Thus, by the definition of Hausdorff dimension to measures, $\rho
_{00}((I\setminus \cup _{j}A_{\theta _{0}^{j}})\cap \cdot )$ has exact
dimension 
\begin{equation}
\alpha =\min_{j}(1-\gamma _{j})  \label{alpdim}
\end{equation}%
which, together with the definition of $\gamma _{j}$, is exactly (\ref{HDAL}%
).

We now replace the Pr\"{u}fer angles $\left( \theta _{n}^{j}\right) _{n\geq
0}$ by a sequence $\left( \zeta _{n}^{j}\right) _{n\geq 0}$ of continuous
piecewise linear functions $\zeta _{n}^{j}=\zeta _{n}^{j}(\varphi _{j})$
which can be shown to be uniformly distributed $\text{mod}$ $\pi $ by the
general metric criterion (see Section 5 of \cite{MarWre}) and whose
difference of their respective Birkhoff average 
\begin{equation*}
E=\frac{1}{N}\sum_{n=1}^{N}\left( f(\theta _{n}^{j})-f(\zeta
_{n}^{j})\right) ~,
\end{equation*}%
for any uniformly continuous function $f$ defined in $\left[ 0,\pi \right] $%
, can be made arbitrarily small by taking the sparseness parameter $\beta $
sufficiently large (see Theorem 5.6 of \cite{MarWre}). As a consequence, (%
\ref{rtr}) is replaced by 
\begin{equation*}
C_{n}^{-1}e^{n(\ln r_{j}-2\left\vert E\right\vert )/2}\leq t_{j,n}\leq
C_{n}e^{n(\ln r_{j}+2\left\vert E\right\vert )/2}
\end{equation*}%
and equations (\ref{c4.4c}) and (\ref{c4.5b}) are affected only by an $%
\varepsilon $ uncertainty, leading to (\ref{epsilon}).

Finally, by the theory of rank one perturbations, we know that $\rho
_{00}(\cup _{j}A_{\theta _{0}^{j}})=0$ holds for almost every $\phi $, and
so for almost every $\phi $ the restriction $\rho _{00}(I\cap \cdot )$ has (%
\ref{HDAL}) as its Hausdorff dimension. This concludes the proof of Theorem %
\ref{thethe}.\newline

\hfill $\Box $


An interesting conclusion drawn from Theorem \ref{thethe} is that the
spectral measure $\rho _{00}$ always inherits the most singular behavior
between its components. Let us explain what this assertion means.

Let $B$ be a Borel set, $B\subset I$ ($I$ given by (\ref{INTI})). If $\alpha
_{\rho _{00}}(\lambda )>0$ for every $\lambda \in B$, then $\rho _{00}(B\cap
\cdot )$ is purely singular-continuous. We see from (\ref{HDAL}) and (\ref%
{rj}) that this holds if, and only if, 
\begin{equation}
(4-\lambda _{j}^{2})(\beta -1)>\left( \frac{1-p^{2}}{p}\right) ^{2}
\label{dbp}
\end{equation}%
\noindent is satisfied for every $j\in \mathcal{I}(\lambda )$. This is
exactly the expression (4.30) of \cite{MarWre}, which gives a necessary
condition for the existence of singular-continuous spectrum (the result
follows directly from Theorem 2.1 of \cite{SS} and Theorem 3.2 of \cite{LS}).

Thus, if condition (\ref{dbp}) fails to be satisfied for at least one $j$ in
some Borel set $B$, then the spectrum of $\rho _{00}(B\cap \cdot )$ if
singular-continuous, it has $0$ Hausdorff dimension. This result is a direct
consequence of Corollary \ref{C4.4a}.


\begin{thebibliography}{DJLS}
\bibitem[DJLS]{DRJLS} R. Del Rio, S. Jitomirskaya, Y. Last and B. Simon,
``Operators with singular continuous spectrum, IV. Hausdorff dimension, rank
one perturbations and localization'', \textit{J. Anal. Math.} \textbf{69},
153-200 (1996).

\bibitem[F]{Falconer} K. J. Falconer, \textquotedblleft Fractal
Geometry\textquotedblright , Wiley, Chichester 1990

\bibitem[GP]{GP} D. J. Gilbert and D. B. Pearson, \textquotedblleft On
subordinacy and analysis of the spectrum of one--dimensional Schr\"{o}dinger
operators\textquotedblright , \textit{J. Math. Anal. Appl.} \textbf{128},
30-56 (1987).

\bibitem[JL]{JitLast} Svetlana Jitomirskaya and Yoram Last,
\textquotedblleft Power--law subordinacy and singular spectra I. Half--line
operators\textquotedblright\ \textit{Acta. Math.} \textbf{183}, 171-189
(1999).

\bibitem[KLS]{KLS} Alexander Kiselev, Yoram Last and Barry Simon. ``Modified
Pr\"{u}fer and EFGP transforms and the spectral analysis of one--dimensional
Schr\"{o}dinger operators'', Commun. Math. Phys. \textbf{194}, 1-45 (1998)

\bibitem[KR]{KR} Denis Krutikov and Christian Remling. \textquotedblleft Schr%
\"{o}dinger Operators with Sparse Potentials: Asymptotics of the Fourier
Transform of the Spectral Measure\textquotedblright , Commun. Math. Phys. 
\textbf{223}, 509-532 (2001)

\bibitem[La]{LANCASTER} Lancaster, P., Tismenetsky, M., \textquotedblleft
The Theory of Matrices\textquotedblright , Academic Press (San Diego),
second edition (1985).

\bibitem[L]{Last} Yoram Last, ``Quantum dynamics and decomposition of
singular continuous spectra'' \textit{J. Funct. Anal.} \textbf{142}, 406-445
(1996).

\bibitem[LS]{LS} Yoram Last and Barry Simon. \textquotedblleft
Eigenfunctions, transfer matrices, and absolutely continuous spectrum of
one--dimensional Schr\"{o}dinger operators\textquotedblright , Invent. Math. 
\textbf{135}, 329-367 (1999)

\bibitem[MWGA]{MarWre} D. H. U. Marchetti, W. F. Wreszinski, L. F. Guidi and
R. M. Angelo, ``Spectral transition in a sparse model and a class of
nonlinear dynamical systems'' \textit{Nonlinearity} \textbf{20}, 765-787
(2007).

\bibitem[P]{Pearson} D. B. Pearson. \textquotedblleft Singular Continuous
Measures in the Scatering Theory\textquotedblright , Commun. Math. Phys. 
\textbf{60}, 13-36 (1978)

\bibitem[P1]{P1} D. B. Pearson, \textquotedblleft Value distribution and
spectral analysis of differential operators\textquotedblright , J. Phys. 
\textbf{A}: Math. Gen. \textbf{26}, 4067-4080 (1993)

\bibitem[RS]{RESI} M. Reed and B. Simon, \textquotedblleft Methods of Modern
Mathematical Physics I: Functional Analysis\textquotedblright , Academic
Press (New York), second edition (1995).

\bibitem[RT]{RodTay} C. A. Rodgers and S. J. Taylor, ``The analysis of
additive set functions in Euclidean space'', \textit{Acta. Math.} \textbf{101%
}, 273-302 (1959).

\bibitem[S]{Saks} S. Saks, \textquotedblleft Theory of the
integral\textquotedblright , Hafner (New York), second edition (1937).

\bibitem[SS]{SS} B. Simon and G. Stolz. \textquotedblleft Operators with
singular continuous spectrum. V. Sparse potentials\textquotedblright , Proc.
Amer. Math. Soc. \textbf{124}, 2073-2080 (1996)

\bibitem[T]{LASEIEU} G. Teschl, \textquotedblleft Jacobi Operators and
Completely Integrable Nonlinear Lattices\textquotedblright , AMS (New York),
(2000).

\bibitem[Z]{Zla} Andrej Zlatlo\v{s}. \textquotedblleft Sparse potentials
with fractional Hausdorff dimension\textquotedblright , J. Funct. Anal. 
\textbf{207}, 216-252 (2004)
\end{thebibliography}
\end{document}